\documentclass[12pt,reqno]{amsart}

\usepackage%[dvips] %including this option causes pdflatex to fail to
{graphicx} 

\usepackage{amssymb,  %this is needed for \measuredhornangle
hyperref
}

\usepackage{breakurl}   %this allows for url breaks in pdf

\usepackage[export]{adjustbox}  %this is necessary for graphics, figures

\usepackage{epigraph}

\usepackage{setspace}  %this is needed for double-spacing

\theoremstyle{definition}

\numberwithin{equation}{section}

\newcommand\N {{\mathbb N}} 

\newcommand\R {{\mathbb R}}

\newcommand\st{{\rm st}}

\title{A philosophical history of infinitesimals}

\author[V. Kanovei]{Vladimir Kanovei} \address{IITP RAS, Moscow,
  Russia}\email{kanovei@googlemail.com}

\author[M. Katz]{Mikhail G. Katz}\address{Department of Mathematics,
  Bar Ilan University, Ramat Gan 5290002
  Israel}\email{katzmik@math.biu.ac.il}

\author[T. Kudryk]{Taras Kudryk} \address{Department of Mathematics,
  Lviv National University, Lviv, Ukraine}\email{kudryk@mail.lviv.ua}

\author[K. Kuhlemann]{Karl Kuhlemann}\address{Gottfried Wilhelm
  Leibniz University Hannover, D-30167 Hannover,
  Germany}\email{kus.kuhlemann@t-online.de}


\begin{document}

%\begin{center}10000 WORDS\end{center}

\begin{abstract}
We explore the issue of providing a foundational framework for
Leibnizian infinitesimals in the light of modern standard and
nonstandard \mbox{approaches}.  We outline a trichotomy of ordinals,
cardinals and ringinals as a historiographic tool.  A \emph{ringinal}
is a concept of infinite number, arithmetic in nature, different from
Cantor's transfinite ordinals and cardinals.  The continuum is not
necessarily identifiable with~$\R$; even if one seeks such an
identification, infinitesimals are not ruled out.  Analysis with
unlimited numbers (via the predicate \emph{standard}) is possible in a
conservative extension of Zermelo--Fraenkel set theory and in this
sense is epistemologically `safe'.  We sketch a recent theory of
infinitesimal analysis that formalizes Leibnizian definitions and
heuristic principles while eschewing both the axiom of choice and
ultrafilters, thus challenging received philosophical views on the
nature of infinitesimals.
%160 words
\end{abstract}

%\doublespacing

\thispagestyle{empty}

%\huge


\maketitle
\tableofcontents

%\today

\epigraph{If in the following someone will complain about the use of
  these quantities, he will show himself to be either an ignorant or
  an ungrateful person.  As an ignorant one, if he does not understand
  what a great light is lit here in the whole method of indivisibles
  and in the field of quadratures; as an ungrateful one, if he hides
  the benefit he gets.\\--Leibniz, \emph{De Quadratura Arithmetica},
  Scholium to Propositio~XXIII}

\section{Intermediate state between Leibniz and d'Alembert}  %1

Abraham Robinson included a sketch of the history of infinitesimal
analysis in Chapter 10 of his 1966 book \cite{Ro66}.  Over half a
century later, significant advances have taken place both in
understanding this history and in axiomatic nonstandard analysis.

Indivisibles were controversial in the 17th century, not least in the
eyes of the jesuits who issued numerous bans against them.%
\footnote{McCue \cite[p.\;419]{Mc68}; Festa \cite[p.\,101]{Fe91} and
\cite{Fe92}.}
%
%and a recent summary in \cite{24c}.  
%
Perhaps sensing the doctrinal burden of the (overly) evocative term
\emph{indivisible} (see Section~\ref{f1}), Leibniz coined the term
\emph{infinitesimal} in 1673 taking up a suggestion of Nicolaus
Mercator.%
\footnote{See Probst \cite[p.\;200]{Pr18}, who also notes that Wallis
  used the term \emph{pars infinitesima} already in 1670.}
Ever since the inception of infinitesimal calculus in the work of
Leibniz and Newton, there have been attempts to delegitimize
infinitesimals.  Abraham Robinson claimed to have developed a
formalisation of Leibnizian infinitesimals.  Robinson's claim has been
challenged by some historians and supported by others.

\subsection{Law of Continuity and \emph{status transitus}}
\label{s13}

Leibniz formulated a number of heuristic principles that governed his
infinitesimal calculus, such as the \emph{Law of \mbox{Continuity}}.%
\footnote{See e.g., Bos \cite{Bo74}, Katz and Sherry \cite{12e}.}
One of the formulations of the Law of Continuity involves the
\mbox{postulation} of a \emph{status transitus}, an intermediate stage
before reaching the end of the process of vanishing.%
\footnote{See Bair et al.~\cite[Section~6.1]{21a}.}
Such an intermediate stage witnesses the \mbox{appearance} of
\emph{inassignable} quantities such as infinitesimals.  The Leibnizian
distinction \mbox{between} assignable and inassignable quantity goes
back to Nicholas of Cusa (Cusanus) (1401--1464).%
\footnote{See Knobloch \cite[p.\;5]{Kn25}.}
A (positive) infinitesimal is smaller than every assignable number
(see further in Section~\ref{s75}).  Concerning such quantities,
Leibniz wrote: 
\begin{enumerate}\item[]
Although they are not assignable, they turn out to be something
existing and not an absolute nothing.%
\footnote{``Bien qu'elles ne soient pas assignables, elles se
  trouvent \^etre quelque chose d’existant et non pas un rien
  absolu.'' (Leibniz as quoted in Bella \cite[p.\,195]{Be19}).}
\end{enumerate}
It is ``a question of fundamental methodology''%
\footnote{Cf.~Archibald et al.~\cite{Ar22b}.}
in the Leibnizian calculus that such infinitesimal quantities are
fictional; see \cite{24c}.  Another formulation of the Law of
Continuity posits that
\begin{enumerate}\item[]
The rules of the finite are found to succeed in the infinite and vice
versa.%
\footnote{``il se trouve que les regles du fini reussissent dans
l'infini, {\ldots} et que vice versa les regles de l'infini
reussissent dans le fini'' Leibniz \cite[pp.\;93--94]{Le02} (original
spelling preserved).  Cf.\;Robinson \cite[p.\;266]{Ro66}.  See further
in Section~\ref{s62}.}
\end{enumerate}
Leibniz's dealing in such notions was not uniformly accepted by his
contemporaries.%
\footnote{Today, we possess a better appreciation of Leibniz's work,
and few scholars would accept \mbox{Whiteside's} claim that ``very few
proofs of any kind in classical mathematics will be allowable, and
certainly none were given in the 17th century on any but the most
elementary numerical level'' \cite[p.\,184]{Wh61}.}
One of the opponents was Michel Rolle, who initiated a lively debate
at the French Academy.%
\footnote{Mancosu \cite{Ma89}; Bair et al.~\cite[Section~2.9]{18a}.}

The 17th century acrimonious debates over indivisibles set the tone
for centuries to come.  In 18th century France, Jean D'Alembert
reasoned as follows:
\begin{enumerate}\item[]
A quantity is either something or nothing: if it is something, it has
not yet vanished; if it is nothing, it has literally vanished.  The
supposition that there is an intermediate state between these two is a
\emph{chimera}.%
\footnote{D'Alembert as translated by Boyer \cite[p.\;248]{Bo59};
  emphasis added.  See Lamand\'e \cite[p.\;52]{La19} for the
  original.}
\end{enumerate}
D'Alembert's attitude toward the intermediate state appears to be less
positive than Leibniz's.  Some historians today adopt the attitude of
the former rather than the latter; see e.g., Section~\ref{s43}.

\subsection{Theological context}
\label{f1}

Indivisibles and atomism were thought contrary to catholic doctrine as
established by the Council of Trent in the 16th century, and
specifically Session 13, canon\;2.  The said canon concludes that
whoever denies the \mbox{transubstantiation} interpretation of the
eucharist, `let him be anathema'.
%
%\footnote{For practical implementations of such an approach, see
%see \cite[Section~4.2]{24c}.  
%
\mbox{Opposition} to atomism within catholicism goes back at least to
Duns Scotus in the 13th \mbox{century}.%
\footnote{Cross \cite[p.\,118]{Cr98}.  and Katz et al.~\cite{23g}.}
Theological concerns may have been behind Fermat's caution in
presenting his method of \emph{adequality} and the variable $E$.%
\footnote{See \cite{13e}, \cite{18d}.}

Leibniz was keenly aware of such doctrinal tensions over indivisibles.
Already in 1668, he published a work on substance aimed
\begin{enumerate}\item[]
to effect a reconciliation between Roman Catholics and Protestants.
{\ldots} These works are especially valuable for what they reveal
about the \mbox{motivations} behind Leibniz's first account of
substance.%
\footnote{Mercer and Sleigh \cite[p.\;68]{Me94}.}
\end{enumerate}
Decades later, in a letter to des Bosses dated 8\;september 1709,
Leibniz \mbox{explicitly} distanced himself from both
transubstantiation and consubstantiation, and sketched a monad-based
approach.%
\footnote{See Look and Rutherford (transl.)~\cite[p.\,153]{Le07}.  On
the theological background, see further in Katz et al.~\cite{24c}.}

\subsection{Ordinals, Cardinals, Ringinals}
\label{s42}

Some historians today tend to interpret the thrust of their training
in naive set theory as entailing that there could exist only two types
of infinite number: \emph{ordinal} and \emph{cardinal}; for examples,
see Sections~\ref{s43} and \ref{s34}.%
\footnote{For a discussion of teleology-prone aspects of modern
historiography of mathematics, see Bair et al.~\cite[Section~4]{22a}.}
Understandably, they are puzzled by the appearance of infinite numbers
in modern frameworks for analysis with infinitesimals (see
Section~\ref{s19}).

Unlike ordinals and cardinals, such infinite numbers are naturally
(nonstandard) elements of the \emph{ring} of integers, or the semiring
of natural numbers.  \mbox{Accordingly}, such infinite numbers could
be referred to as \emph{ringinals},%
\footnote{Note that \emph{ringinal} is a different concept from
\emph{numerosity}.  The existence of a numerosity function is
%equivalent to the existence of a selective ultrafilter, hence it is
independent of the axioms of ZFC \cite[Abstract]{Be03}, whereas the
existence of ringinals is provable in ZFC and is, moreover,
conservative over ZF~\cite{21e}.}
to emphasize the \mbox{contrast} with Cantorian infinities.  Such a
ringinal is necessarily greater than every naive counting
number~$1,2,3,\ldots$\, Laugwitz observed that Cantor had opposed such
\mbox{entities}:
\begin{enumerate}\item[]
%Being 
[C]onvinced that he had discovered the one and only way to establish
the actual infinite in mathematics after millenia of philosophical
prejudices, Cantor was irritated by attempts of others to grasp the
infinitely large, or to revive the infinitely small.%
\footnote{Laugwitz \cite[p.\,102]{La02}.}
\end{enumerate}
Our distinction between Cantorian infinities and ringinals is related
to Laugwitz's distinction ``between Cantor's transfinite arithmetic
and the theory of ordered \mbox{algebraic} structures.''%
\footnote{Ibid.}
As noted by Laugwitz, 
\begin{enumerate}\item[]
To Stolz, Veronese and Levi-Civita we owe early insights in ordered
\mbox{mathematical} structures as a \emph{third} aspect of the number
concept.%
\footnote{Laugwitz \cite[p.\,103]{La02}; emphasis added.}
\end{enumerate}
The three aspects alluded to by Laugwitz are, in our terminology,
ordinal, cardinal, and ringinal.

\subsection{Leibniz's bounded infinity and Knobloch's cardinals}
\label{s43}

With regard to Cantor's denial of the possibility of ringinals (see
Section~\ref{s42}), some modern \mbox{historians} appear to follow
suit.  As a case study, we analyze Knobloch's comments on
\mbox{Leibnizian} infinitesimal and infinite numbers, made in the
context of an analysis of \emph{Propositio XI} of \emph{De Quadratura
Arithmetica} (see e.g., \cite{Le04b}).  In the Leibnizian
passage,~$\lambda$ denotes the curve admitting the axis of the
ordinates as an asymptote,%
\footnote{In the context of \emph{Propositio XI}, the curve is of type
$y=\frac1{x^2+1}$ or more precisely~$x=\frac1{y^2+1}$.  In Leibniz's
figures, the axis of the ordinates is horizontal.}
whereas~$(\mu)$ denotes an infinitesimal abscissa, and~$(\mu)\lambda$
denotes the corresponding ordinate which is necessarily infinite.%
\footnote{I.e., an \emph{infinitum terminatum}, or in modern terms
\emph{unlimited}; see Section~\ref{s52}.}
In this connection, Knobloch writes:
\begin{enumerate}
\item[[1\!\!]]  Leibniz uses also another terminology:
  \emph{ordinata~$(\mu)\lambda$ erit longitudine \mbox{infinita},
  major qualibet assignabili~${}_4B\,_4D$} (the ordinate will be of
  infinite length, larger than any assignable ordinate~${}_4B\,_4D$)%
\footnote{We use \emph{subscripts} on the left side, as in~${}_4B$, to
make the formulas easier to read for the modern reader.  In his
manuscripts, Leibniz placed them on the same line as the letter, as in
{\tiny4}$B$.}
{\ldots}
\item[[2\!\!]]  Hence there is an unavoidable consequence.  The set of
  all finite cardinal numbers~$1, 2,3,\ldots$~is a \emph{transfinite
  set}.  Its cardinal number is Alef$_0$.%
\footnote{Here Alef$_0$ is Knobloch's notation.}
This is the least cardinal number being [sic] larger than any finite
cardinal number.
\item[[3\!\!]]  Leibniz's terminology implies actual infinity though
  he rejects the \mbox{existence} of an infinite number, etc.%
\footnote{Knobloch \cite[pp.\;13--14]{Kn18}; emphasis on ``transfinite
  set'' and numerals [1], [2], [3] added.}
\end{enumerate}
Knobloch's analysis starts (part [1]) and ends (part [3]) with
Leibniz, suggesting that the intermediate discussion (part [2]) of
``transfinite sets'' and Alef$_0$ is meant to shed light on Leibniz's
mathematics.  Yet the reader may well wonder what transfinite sets and
the Cantorian Alef$_0$ have to do with Leibniz's mathematics.

Note the symbol~$(\mu)\lambda$ used in Knobloch's passage. Leibniz
exploited this symbol to denote a bounded infinity (\emph{infinitum
terminatum}); see further in Section~\ref{s73}.  It is instructive to
analyze Knobloch's reasoning in this pasage.  Knobloch \mbox{appears}
to \mbox{argue} that since Leibniz's bounded infinity~$(\mu)\lambda$
is a magnitude greater than the naive integers~$1,2,3,\ldots$, the
`unavoidable consequence' would be that it must be at least Alef$_0$.
But is such a consequence unavoidable?  Besides the fact that
Knobloch's argumentation would be meaningless to Leibniz who holds
infinite wholes (such as the Cantorian Alef$_0$) to be contradictory
because contrary to the part-whole principle, Knobloch's argument
depends on an implicit assumption that any magnitude would necessarily
be below some naive counting number~$1,2,3,\ldots$\, In effect,
Knobloch is denying the possibility of ringinals.

Knobloch claims that ``Leibniz's terminology implies \emph{actual
infinity} though he rejects the existence of an infinite number''%
\footnote{Emphasis added.}
in the passage above.  But which ``actual infinity'' is Knobloch
referring to exactly?  Here Knobloch similarly fails to distinguish
between infinite magnitude (namely, the Leibnizian \emph{infinita
terminata}) and infinite multitude (namely, infinite wholes which
Leibniz indeed rejected as contradictory).

Knobloch's interpretation of Leibniz ultimately leads him to attribute
\mbox{contraditions} to Leibniz where there may be none.  As analyzed
in \cite{23h}, there is no contradiction between Theorems 11 and 45 in
Leibniz's \emph{De Quadratura Arithmetica}, as they deal with
different notions of infinity: Theorem 11 deals with the
\emph{infinitum terminatum}, whereas Theorem 45 deals with an ideal
perspective point at infinity.  In a recent text, Knobloch discusses
both theorems and concludes as follows:
\begin{enumerate}\item[]
Here [i.e., in Theorem~45], Leibniz says the opposite of what he had
said in the demonstration of theorem 11.  There, he had explicitly
excluded the possibility that the curve meets the asymptote in
accordance with the meaning of its name.%
\footnote{Knobloch \cite[p.\;9]{Kn24}.}
\end{enumerate}
It emerges that Knobloch prefers to attribute a contradiction to
Leibniz rather than accept an interpretation differing from Ishiguro's
as developed in \cite[Chapter\;5]{Is90}.%
\footnote{Knobloch adds the following sentence: ``We shall come back
  to this matter of fact'' (ibid.), but he never does.}

Meanwhile, Leibniz made it clear that his magnitudes violated the
concept of comparability expressed in Euclid V, Definition 4, closely
related to the Archimedean principle.  This enabled Leibniz to use
bounded infinities while adhering to the part-whole principle; for
details see Section~\ref{s75}.

\subsection{Sonar on hornangles}
\label{s34}

\begin{figure}\begin{center}\includegraphics [scale=0.53,rotate=0,trim = 0 0 0 0,clip]{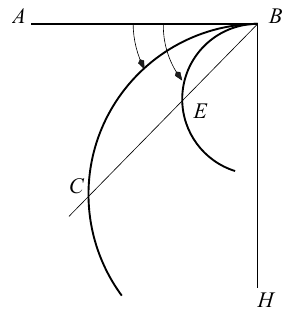}\end{center}\caption{Sonar's analysis of hornangles}
\label{figure1}
\end{figure}

Related misconceptions about non-Archi\-medean \mbox{behavior} occur
in Sonar's work.  Sonar analyzes de Saint Vincent's discussion of
hornangles denoted~$\measuredangle ABC$ and~$\measuredangle ABE$.%
\footnote{Sonar \cite[p.\;227]{So21}.}
With reference to Figure~\ref{figure1}, he shows that the hornangles
must be necessarily smaller than a rectilinear angle~$\measuredangle
ABH$ as well as smaller than its fractions:
\begin{equation}
\label{e41}
\frac{\measuredangle ABH}2 > \frac{\measuredangle ABH}{2^2} > \ldots >
\frac{\measuredangle ABH}{2^n} > \ldots > \measuredangle ABC
\end{equation}
However, the conclusion Sonar draws at this point is revealing of his
assumptions, as he writes on page 228:
\begin{quote}
But if~$n$ increases~$\measuredangle ABH/2^n$ decreases and
\emph{hence it becomes clear} that
\begin{equation}
\label{e42}
\measuredangle ABC = \measuredangle ABE = 0
\end{equation}
holds for the horn angles.  It is rather fascinating that Gr\'egoire
de Saint-Vincent realised this fact but that he could never accept it!%
\footnote{Sonar \cite[p.\;228]{So21}; emphasis added.}
\end{quote}
On page 228, Sonar apparently considers himself to have proved the
claim that hornangles are zero, refers to such a claim as a `fact',
and comments that de Saint Vincent never accepted such a fact.  Note
that de Saint Vincent's teacher Clavius was favorably disposed toward
hornangles and was involved in a controversy with Peletier on this
topic.%
\footnote{Maier\`u \cite{Ma90}, Malet \cite{Ma97}, Axworthy
  \cite[pp.\,12--13]{Ax18} and references therein.}
As acknowledged by Sonar,
\begin{enumerate}\item[]
Gottfried Wilhelm Leibniz will resume de Saint-Vincent’s work later.%
\footnote{Sonar \cite[p.\;228]{So21}.}
\end{enumerate}
Indeed, Leibniz also used hornangles.%
\footnote{For a discussion of hornangles in Leibniz, see e.g.,
\cite[Section~4.5]{24c}.}
As late as 1826, Cauchy spoke of hornangles (\emph{angles de
contingence}) in his work on differential geometry \cite{Ca26}; see
\cite{17a}, \cite{18e}, \cite{19a}, \cite{21f},
\cite[Section~3.10]{22a}.
%Section~\ref{s31b}.

Thus, Sonar's argument for the vanishing of the hornangle seems to be
based on a hidden hypothesis.  Indeed, the inference from
inequality~\eqref{e41} to identity~\eqref{e42} (described by Sonar's
phrase ``hence it becomes clear that'') implicitly relies upon a
property equivalent to the Archimedean axiom: if the
hornangle~$\measuredangle ABC$ \mbox{satisfies} the bound
$\measuredangle ABC<\frac1n$ for each naive integer~$n$,
then~$\measuredangle ABC=0$.  What is therefore `clear' is that Sonar
is relying on the Archimedean axiom.  If one refrains from relying
upon the Archimedean axiom, there is no reason to assume that
hornangles vanish.  Sonar's inference is another instance of a denial
of the possibility of ringinals (cf. Section~\ref{s43}).

Earlier in his book, Sonar writes on page 40:
\begin{enumerate}\item[]
[A]lready Eudoxus knew that also other number systems -- so-called
non-Archimedean number systems -- were conceivable {\ldots} Such a
system of quantities which was already known to the Greeks were
cornicular angles or horn angles.%
\footnote{Sonar \cite[p.\;40]{So21}.}
\end{enumerate}
Accordingly, Sonar is prepared to envision hornangles as part of a
non-Archimedean number system.  It is unclear how his comment on
page\;40 is to be squared with his comments on page 228 where
hornangles are declared to provably vanish.  Sonar's claim that ``de
Saint-Vincent realized [that hornangles vanish] but could never
\mbox{accept} it'' must be rejected as a presentist interpretation.

We examine an additional case of presentist treatment of historical
\mbox{infinitesimalists} in Section~\ref{s15}.

\subsection{Cauchy's infinitesimal delta function}
\label{s15}

In 1982, L\"utzen discusses delta functions that occurred in Poisson
and Cauchy in connection with Fourier's integral theorem.
Specifically, he mentions what is called today the Cauchy probability
density.  L\"utzen writes:
\begin{enumerate}\item[]
They both used arguments which in distribution language would be
\[
\frac12 \, \frac{\alpha}{\alpha^2+(x-a)^2} \;\to\; \delta(x-a), \quad
\text{for } \alpha\to 0.%
\footnote{L\"utzen \cite[p.\,115]{Lu82}.}
\]
\end{enumerate}
Note that the meaning of the first arrow would require explanation `in
distribution language', which Cauchy did not possess.  L\"utzen's
formulation is rather presentist, since he does not even mention
that~$\alpha$ was actually infinitesimal for Cauchy.

As early as 1971, Hans Freudenthal briefly mentioned Cauchy's work on
\begin{quote}
\mbox{singular} integrals (i.e., integrals of infinitely large
functions over infinitely small paths [$\delta$\;functions]).%
\footnote{Freudenthal \cite [p.\,135] {Fr71}.}
\end{quote}

A 1989 article by Detlef Laugwitz contains a section 5.5 starting on
page 227 entitled ``Cauchy and delta functions".  Laugwitz mentions
Cauchy's use of the ``\mbox{language} of infinitesimals''%
\footnote{Laugwitz \cite[p.\;229]{La89}.}
on page 289 of Note XVIII in Cauchy's \mbox{publication} from 1827
\cite{Ca27}.  This is in reference to Cauchy's infinitesimal~$\alpha$
appearing in the \mbox{expression}
\[
\frac{\alpha}{\alpha^2+(x-a)^2}
\]
which integrated against~$F(x)$ produces the value~$F(a)$ of~$F$ at
the point~$a$ (up to a factor),%
\footnote{From the modern point of view, the relation is not that of
equality but rather infinite proximity.  For further details see Katz
and Tall \cite{13g}.}
a property sometimes referred to as the sifting property (of the delta
function).  More specifically, Cauchy wrote:
\begin{quote}
Moreover one finds, denoting by~$\alpha$,~$\epsilon$ two infinitely
small numbers,
\[
\frac{1}{2} \int_{a-\epsilon}^{a+\epsilon} F(\mu) \frac{\alpha \;
  d\mu}{\alpha^2 + (\mu-a)^2} = \frac{\pi}{2} F(a).%
\footnote{Cauchy \cite[p.\,289]{Ca27}.}
\]
\end{quote}

L\"utzen is by no means the only author who acknowledges Cauchy's role
in the prehistory of distributions while sweeping under the rug his
use of infinitesimals in writing down his delta function.  As early as
1955, van der Pol and Bremmer discuss the prehistory of distributions
and mention Cauchy's role in the history of the delta function
possessing what they refer to as the \emph{sifting property},%
\footnote{van der Pol and Bremmer \cite[p.\;64]{va55}.}
but omit to mention Cauchy's use of infinitesimals in his delta
function.%
\footnote{van der Pol and Bremmer do quote Hermite on infinitesimals
\cite[pp.\;62--63]{va55} and Poisson on infinitesimals
\cite[pp.\;63--64]{va55}.}

\section{Lagging formalisation} 
\label{s18}

The formalisation and/or axiomatisation of analysis achieved during
the decades around 1900 did not incorporate a formalisation of
infinitesimals.  The Zermelo--Fraenkel set theory (ZF) is a theory in
the~$\in$-language (thus,~$\varnothing\in \{\varnothing\}$, etc.).
The theory ZF was completely articulated by the 1920s.  Meanwhile, an
applicable formalisation of infinitesimals lagged by half a century
(as we know in retrospect; see Section~\ref{s19}).  Felix Klein sensed
such a lag, and lamented it in the following terms:
\begin{enumerate}\item[]
The question naturally arises whether \ldots{}\;it would be possible
to \emph{modify} the traditional foundations of infinitesimal
calculus, so as to include actually \emph{infinitely small} quantities
in a way that would satisfy modern demands as to rigor; in other
words, to construct a non-Archime\-dean system.  The first and chief
problem of this analysis would be to prove the mean-value theorem
$f(x+h)-f(x)=h \cdot f'(x+\vartheta h)$ [where~$0\leq\vartheta\leq1$]
from the assumed axioms.  I will not say that progress in this
direction is impossible, but it is true that none of the investigators
have achieved anything positive.%
\footnote{Klein \cite[p.\;219]{Kl08}; emphasis added.  See also
\cite[Section\;6.1]{13c}, \cite{18i}.  Klein has not been treated
fairly by all modern historians; see \cite{18b} and \cite{23b}.}
\end{enumerate}
Progress in this direction started with Skolem \cite{Sk33}; see
Section~\ref{s51c}.

\subsection{Campaign, from Russell to Carnap}
\label{s51b}

In the meantime, a campaign of \mbox{demonisation} of infinitesimals
took place, that was a thinly veiled attempt to \mbox{conceal} the
failure of existing formalisations to incorporate infinitesimals; see
\cite{13a}.

As detailed by Ehrlich in \cite{Eh06}, Bertrand Russell and others
publish texts claiming to derive the inconsistency of infinitesimals
from philosophical principles, which merely embodied the philosophical
prejudices of their authors.%
\footnote{A more sympathetic attitude toward infinitesimals at the
  time was displayed by the neo-Kantian school led by Hermann Cohen;
  see Laugwitz \cite[p.\,111]{La02}; Mormann and Katz \cite{13h};
  Edgar~\cite{Ed20}; Pringe \cite{Pr23}.}
According to Russell,
\begin{enumerate}\item[]
[I]nfinitesimals as explaining continuity must be regarded as
unnecessary, erroneous, and self-contradictory.%
\footnote{Russell \cite[item~324]{Ru03}.}
\end{enumerate}

It is hard to miss the tone of smug satisfaction in Florian Cajori's
1917 remarks concerning
\begin{quote}
\emph{wonderful} strides in the banishment of infinitely small
quantities.%
\footnote{Cajori \cite[p.\,152]{Ca17}; emphasis added.}
\end{quote}
Significantly, Cajori is of the (unsubstantiated) opinion that
``Leibniz's philosophy of the calculus was poor.''%
\footnote{\label{f28}Cajori \cite[p.\,153]{Ca17}.  Cajori notes
further that ``[The] return to the use of infinitely small quantities
is noticeable in several English texts of the second half of the
[18th] century.  An old lady once defended Calvinism by saying that if
you took away her total depravity you took away her religion.  There
were mathematicians who believed that if you took away infinitely
small quantities you took away all their mathematics''
\cite[p.\,151]{Ca17}.  Cajori's suggestion of parallelism between a
lady's `depravity' and infinitesimals illustrates the attitudes toward
infinitesimals prevalent during the 1910s.}
%Recall that Lamarle had already referred to infinitesimals as a
%``vice radical'' (see Section~\ref{s33}).

Rudolf Carnap highlights the purported force of one of his own
philosophical innovations -- the notion of a \emph{pseudo concept} --
as allegedly diagnosing the concept of the infinitesimal.  To clarify
further, Carnap includes a parenthetical epithet \emph{empty words}%
\footnote{Carnap \cite[pp.\;306--307]{Car28}.}
for what Leibniz thought were \emph{useful fictions}.%
\footnote{For details on useful fictions see Sherry and Katz
  \cite{14c}; Bair et al.~\cite[Section~1.7]{21a}; Katz et
  al.~\cite{21g}.}
% \cite{22c}.
%

Writing in \emph{The Monist} in 1925 and possibly taking cue from
Cantor's \mbox{epidemiological} metaphor of `cholera bacillus of
mathematics',%
\footnote{See \cite[p.\;233; note 55, p.\;349]{Da90}.}  Parkhurst and
Kingsland speak of mathematics and metaphysics as territory
``contaminated''%
\footnote{Parkhurst and Kingsland \cite[p.\;633]{Pa25}.}
and ``infected''%
\footnote{Parkhurst and Kingsland \cite[p.\;634]{Pa25}.}
by infinitesimals, which risk causing ``peripatetic fever.''%
\footnote{Ibid.}
Note that Cantor was writing in the midst of the fifth cholera
outbreak (1881--1896), while Parkhurst and Kingsland, on the heels of
the Spanish flu (1918--1920).  The common denominator is the implied
view of infinitesimals as a plague.

\subsection{Courant's duty}
\label{s51}

Richard Courant elevates the \emph{avoidance} of foggy, hazy
\mbox{infinitesimals} to the status of a \emph{duty}:
\begin{enumerate}\item[]
[W]e must beware of regarding the derivative as the quotient of two
\mbox{quantities} which are actually ``infinitely small''. {\ldots} It
is true that this fog that hung round the foundations of the new
science did not prevent Leibnitz or his great successors from finding
the right path.  But this does not release us from \emph{the\;duty of
avoiding} every such hazy idea in our building-up of the differential
and integral calculus.%
\footnote{Courant \cite[p.\,101]{Co37}; emphasis added.}
\end{enumerate}
Courant does not attempt to explain how exactly Leibniz and `his great
successors' managed to ``find the right path'' while dealing in foggy,
hazy entities (but see Section~\ref{s53b}).  Courant's admonition is
assorted with another \emph{must} and is rather specific:
\begin{enumerate}\item[]
We must {\ldots}\;guard ourselves against thinking of~$dx$ as an
``infinitely small quantity'' or ``infinitesimal'', or of the integral
as the ``sum of an infinite number of infinitely small quantities.''%
\footnote{Courant \cite[pp.\;80--81]{Co37}.}
\end{enumerate}

Echoing Carnap's \emph{empty words} (see Section~\ref{s51b}), in 2021
Costantini devotes his Section~7 to an analysis of ``The case of
infinitesimals and other \emph{empty~notions}.''%
\footnote{Costantini \cite[pp.\;285--286]{Cos21}; emphasis added.}

It apparently occurred to few of these 20th and 21st century scholars
that, \mbox{arguably}, what we \emph{must} develop is an appropriate
formalisation of the Leibnizian distinction between assignable and
inassignable quantities (see Sections~\ref{s13} and~\ref{s62}).

\subsection{Finding the right path}
\label{s53b}

How \emph{did} the founders of the calculus manage to ``find the right
path'' as Courant put it (see Section~\ref{s51})?  Perhaps the
entities in question were not as hazy as Courant made them out to be.
Courant's original comments date from 1927 \cite{Co27}.  Nearly half a
century later, historian Bos will suggest that
\begin{enumerate}\item[]
A~preliminary explanation of why the calculus could develop on the
\mbox{insecure} foundation of the acceptance of infinitely small and
infinitely large \mbox{quantities} is provided by the recently
developed \emph{non-standard analysis}, which shows that it is
possible to remove the inconsistencies without removing the
\mbox{infinitesimals} themselves,%
\footnote{Bos \cite[p.\,13]{Bo74}.}
\end{enumerate}
but seems to walk back his `preliminary explanation' seventy pages
later.%
\footnote{Bos \cite[pp.\;81--83]{Bo74}.  For an analysis of Bos's
  remarks, see Katz and Sherry \cite[Section\;11.3]{13f} and Bair et
  al.~\cite[Section~2.7]{17b}.}

\subsection{Of darts and chimeras}
\label{s53}

In 1969, Bernstein and Wattenberg published a paper applying
Robinson's framework to measure theory.  The paper opens with a
discussion of the dart experiment (probabilistic analysis of throwing
a dart at a target).%
\footnote{Bernstein and Wattenberg \cite[p.\,171]{Be69}.}
The dart paradox consists in the observation that every point of the
target has zero chance of being hit, yet some point does get hit by
the dart.  Bernstein and Wattenberg argued that Robinson's
infinitesimal analysis enables a resolution of the paradox.  They
exploited a hyperfinite set containing~$\R$; the resulting weighted
counting measure assigns a nonzero infinitesimal probability to each
real number.%
\footnote{\label{f91}In axiomatic approaches to infinitesimal
  analysis, one exploits a finite set containing all standard real
  numbers.  The existence of such a set is proved by the theory
  BSPT$'$, which is a subtheory of IST and HST.  The theory BSPT$'$ is
  conservative over ZF, as proved by Hrbacek and
  Katz~\cite[p.\,10]{21e}.  See further in Section~\ref{s52}.}
Bernstein and Wattenberg presented the dart experiment as a feature of
Robinson's approach.

In 1970, Alain Connes publishes a paper on ultrapowers and nonstandard
\mbox{analysis} \cite{Co70} which cites Bernstein and Wattenberg (and
even claims to improve on some of their results), indicating that
Connes was familiar with their treatment of the dart paradox,
\emph{and appreciated it}.

By 2000, Connes seeks to present the dart experiment as a shortcoming
of \mbox{Robinson's} framework.%
\footnote{Connes \cite[pp.\,13--14]{Co00}.}
%See the quote in Section~\ref{s2}.
%
The critical comments by Connes were likely in response to Bernstein
and Wattenberg.

Two decades later in 2021, Connes is still panning ultrafilters.  He
expresses the following two sentiments:
\begin{enumerate}
\item
\label{i1}
a criticism of nonprincipal ultrafilters, and 
\item
an endorsement of the Continuum Hypothesis as a tool for producing a
Dixmier trace with desirable properties.%
\footnote{Connes \cite[p.\;47]{Co21}.}
\end{enumerate}
As a source for the construction of the desirable trace (using the
``medial limit of Mokobodzki''),%
\footnote{Note that Larson \cite{La09} showed in 2009 that ZFC does
  not prove the existence of medial limits.}
Connes cites Meyer \cite{Me73}.  Here Meyer wrote:
%Paul-André Meyer (21 August 1934 – 30 January 2003) 
\begin{enumerate}\item[]
Soit~$E$ l'ensemble de tous les ordinaux d\'enombrables:~$E$ a la
puissance du continu (axiom du continu!), etc.%
\footnote{Meyer \cite[p.\;200]{Me73}.}
\end{enumerate}
Note that this form of the Continuum Hypothesis implies, over ZF (see
the \mbox{beginning} of Section~\ref{s18}), the existence of a
well-ordering of the reals.  This foundational \mbox{material} is
exploited by Connes's source Meyer in the construction of the
desirable trace.  But such a well-ordering suffices to prove the
existence of a nonprincipal ultrafilter, criticized in
item~\eqref{i1}.  Therefore the position as expressed by Connes in
2021 is mathematically incoherent.%
\footnote{See further in Katz and Leichtnam \cite{13d}, Kanovei et
al.~\cite{15a}, \cite[Section~8.6]{21e}, and Sanders \cite{18l}.  A
related critique of Robinson's framework by Paul Halmos is analyzed in
B{\l}aszczyk et al.~\cite{16b}.}
Furthermore, it turned out recently that nonstandard analysis can be
practiced without ultrafilters.%
\footnote{See Section~\ref{s52} on axiomatic frameworks such as SPOT.}

\subsection{Emotionally charged}

The industrialisation of the calculus programs \mbox{starting} in the
1960s produced a series of textbooks setting a certain classroom
style.  \mbox{Advocates} of infinitesimals argue that such classrooms
create an emotionally charged atmosphere (foreshadowed by the tone of
Courant's remarks quoted in Section~\ref{s51}):

\begin{enumerate}\item[]
Students are often given emotionally charged instructions to avoid
thinking of~$dy/dx$ as a quotient and to conceptualize it as a limit,
even though the formulae of the calculus visibly seem to operate as if
it is a quotient involving symbols that can be shifted around to
change differential equations into integrals.%
\footnote{Tall \cite[p.\;336]{Ta13}.}
\end{enumerate}
In such classrooms, instructors sometimes make deprecating comments at
the \mbox{expense} of infinitesimals, inculcating in the student a
faith that there is something \mbox{fundamentally} wrong with them.
Students also tend to internalize a literal interpretation of the
epithet \emph{real} in the expression \emph{real number}.%
\footnote{Studies of infinitesimal-based pedagogy include Katz and
Polev \cite{17h}, Kuhlemann~\cite{Ku22}, \cite{25d}, and Agnesi
\cite{Ag23}.}

\section{From Fine Hall to modern infinitesimal analysis} 
\label{s19}

Abraham Robinson (1918--1974) was an applied mathematician,%
\footnote{It is to be noted that Abraham Robinson's \emph{Wing Theory}
  \cite{Ro56} is not a branch of model theory, but rather ``an
  admirable compendium of the mathematical theories of the
  aerodynamics of aerofoils and wings'' \cite{Li56}.}
logician, and inventor of nonstandard analysis (NSA).\, Robinson named
his theory
\begin{enumerate}\item[]
\emph{Non-standard Analysis}[,] since it involves and was, in part,
inspired by the so-called Non-standard models of Arithmetic whose
existence was first pointed out by T. Skolem.%
\footnote{Robinson \cite[p.\,vii]{Ro66}; emphasis added.}
\end{enumerate}
Thus the name of the field was influenced by Skolem's work.%
\footnote{It is therefore historically inaccurate to claim, as Loeb
and Wolff do, that ``Robinson chose the name `nonstandard analysis'
because the nonstandard world is used to analyze the standard one''
\cite[p.\,vii]{Loeb2}.  As the passage we quoted from Robinson's book
indicates, his choice of the name `nonstandard analysis' was inspired
by Skolem's models, not by the reason claimed by Loeb and Wolff, who
go on to lodge the following philosophical claim: ``Internal set
theory, on the other hand, works with only the nonstandard world, but
recognizes some elements of that world as being `standard'{}''
(ibid.).  Loeb and Wolff use such a philosophical claim to buttress
their contention concerning the alleged impossibility of carrying out
certain constructions (nonstandard hulls, Loeb measures, etc.)~in
axiomatic approaches to NSA.  They conclude: ``These constructions do
not make sense in internal set theory because there is no standard
world'' (ibid.).  Their contention is incorrect, as shown in Hrbacek
and Katz \cite{23c}, where both nonstandard hulls and Loeb measures
are handled within an axiomatic approach to NSA.  On axiomatic
approaches to NSA, see further in Section~\ref{s52}.}

\subsection{Ehrlich and Dauben on NSA}
\label{s51c}

As noted by Ehrlich,
\begin{enumerate}\item[]
After mathematicians had been taught for decades that a consistent
theory of the calculus based on infinitesimals was impossible, Abraham
Robinson was certainly swimming against the tide when he proved
otherwise.  \cite{Eh22}
\end{enumerate}
Building upon earlier breakthroughs by T. Skolem \cite{Sk33},
E. Hewitt~\cite{He48}, J.\;{\L}o\'s~\cite{Lo55}, and others, Robinson
first introduced NSA in a 1961 article~\cite{Ro61},%
\footnote{Robinson reported that the idea of infinitesimal analysis
  came to him as he was walking into Fine Hall in the fall of 1960.
  This is recounted in the biography of Robinson by Dauben
  \cite[p.\;281]{Da95}.}
and then in the 1966 book \cite{Ro66}.%
\footnote{For an analysis of the relation between the infinitesimals
  of Leibniz and Robinson, see Bair et al.~\cite{21a}.}
Dauben notes that
\begin{enumerate}\item[]
Robinson succeeded in showing the reasonableness of ``redrawing'' the
early history of the calculus to reinstate past views that, cast in
the light of nonstandard analysis, could be seen more clearly.%
\footnote{Dauben \cite[p.\;327]{Da21}.}
\end{enumerate}
Robinson's philosophical position of Formalism is clarified in
\cite{26c}.

As anticipated apprehensively by Courant (see Section~\ref{s51}),~$dx$
is \mbox{infinitesimal},%
\footnote{This is the case in Robinson's framework.  Similar notation
$dx$ is used in traditional non-infinitesimal analysis, where this is
not the case.}
and the definite integral is defined via the sum of
a nonstandard number of \mbox{infinitely} small quantities,%
\footnote{Namely, the integral is the standard number infinitely close
to such a sum.  For further details, see e.g., the textbook by Vakil
\cite{Va11}.}
challenging Russell's claim that ``the so-called \mbox{infinitesimal}
\mbox{calculus} \ldots\ has nothing to do with the infinitesimal.''%
\footnote{Russell \cite[item 308]{Ru03}.}

\subsection{Axiomatisations: IST, HST, SPOT; conservativity}
\label{s52}

In the mid-1970s, Edward Nelson \cite{Ne77} developed an axiomatic, or
syntactic, framework for NSA, called Internal Set Theory (IST).\, At
the same time, Karel Hrbacek \cite{Hr78} developed a different
axiomatic/syn\-tactic framework, now called HST; see the monograph by
Kanovei and Reeken~\cite{Ka04}.

The theory SPOT \cite{21e} (an acronym of its axioms) is a subtheory
of both IST and HST.  

Whereas IST and HST are conservative extensions of ZFC (ZF plus the
\mbox{Axiom} of Choice), the theory SPOT is a conservative extension
of ZF itself.%
\footnote{This feature of SPOT enables a reverse-mathematical analysis
  of the role of the axiom of choice in mathematics with
  infinitesimals, something that was not possible over the classical
  axiomatic theories such as IST and HST that themselves rely on the
  axiom of choice.}
It is a \mbox{theory} in the \st-$\in$-language, where ``\st'' is a
one-place predicate (thus, \st$(x)$ reads ``$x$ is \mbox{standard}'')
formalizing the Leibnizian distinction between assignable and
inassignable quantities.%
\footnote{\label{ff}Note that in axiomatic set theories such as IST,
  HST and SPOT, the standard (or ``naive" as Reeb called them)
  integers do not exhaust the set~$\N$.  Standard integers do not form
  a set.  The separation axiom schema of ZF still holds in SPOT but it
  only applies to~$\in$-formulas, not to \st-$\in$-formulas involving
  the predicate \emph{standard}.  See further in Section~\ref{s54}.}

Since SPOT incorporates the axioms of ZF, the natural numbers~$\N$ and
the real numbers~$\R$ are developed as usual.  A number is
\emph{limited} if its absolute value is smaller than some standard
real number; otherwise it is \emph{unlimited}.%
\footnote{An unlimited element of~$\N$ is what was referred to as a
\emph{ringinal} in Section~\ref{s42}.}
An \mbox{infinitesimal} is a number smaller in absolute value than
every positive standard number.%
\footnote{This definition of an infinitesimal is valid both in
Robinson's original framework and in \mbox{syntactic}/axiomatic
frameworks such as IST, HST, and SPOT.\, It formalizes Leibniz's
definition of infinitesimal as smaller than every assignable number;
see Section~\ref{s13}.}
Two numbers are called \emph{infinitely close} if their difference is
an infinitesimal.

In addition to the axioms of ZF (which apply to all~$\in$-formulas),%
\footnote{As mentioned in note~\ref{ff}, the axiom schema of
separation (comprehension) applies to all \mbox{$\in$-formulas}, but
not to formulas involving the new predicate \st, so that in particular
the standard integers don't form a set.  This feature has been used to
provide an account of the sorites paradox in Dinis and van den Berg
\cite[p.\;255]{Di19}; see also Dinis \cite{Di23}.}
SPOT has the following three axioms.  Here~~$\exists^{st}$ and
$\forall^{st}$ denote quantification over standard entities only,
and~$\N$ is the set of natural numbers as defined in ZF (see above).

\medskip\noindent \textbf{T} (Transfer) Let~$\phi$ be an~$\in$-formula
with standard parameters.  \\ Then~$\forall^{st} x \, \phi(x) \iff
\forall x \, \phi(x)$.

\medskip\noindent \textbf{O} (Nontriviality)~$\exists\nu\in\N \;
\forall^{st} n \in\N \, (n \neq \nu)$.

\medskip\noindent \textbf{SP} (Standard Part) Every limited real is
infinitely close to a unique standard real, called its \emph{shadow}.

\medskip\noindent 
Equivalently,

\medskip\noindent \textbf{SP$'$} \quad~$\forall A\subseteq\N\;
\exists^{st}B\subseteq\N \;\forall^{st}n\in\N\,(n\in B\iff n\in A)$.%
\footnote{Intuitively, the equivalence of SP and SP$'$ can be
  understood in terms of the binary expansion of a real number~$x$.
  If~$A$ denotes the set of ranks at which the corresponding digit
  of~$x$ is~$1$, and~$B$ denotes the set of ranks at which the
  corresponding digit of its shadow sh$(x)$ is~$1$, then the fact
  that~$x$ and sh$(x)$ are infinitely close corresponds to the fact
  that~$A$ and~$B$ agree at all standard (limited) ranks.  The
  detailed argument is slightly more technical due to the
  non-uniqueness of binary expansion; see \cite[Lemma~2.4]{21e}.}

\medskip
The Transfer axiom (schema) is a formalisation of the Leibnizian
\emph{Law of \mbox{Continuity}} (see Sections~\ref{s13} and
\ref{s62}).  For a more detailed introduction to SPOT and related
theories such as SCOT and BSPT$'$ see recent work by Hrbacek and Katz
in \cite{24a}, \cite{21e}, \cite{23d}, and \cite{23e}.

\subsection{Leibniz's laws compared to the axioms of SPOT}
\label{s62}

The simplicity of the additional axioms of SPOT (compared to some of
the more abstruse axioms of ZF) suggests that a formalisation of
infinitesimals may not have been beyond the intellectual prowess of
the pioneers of formalisation from over a century ago, which calls to
mind the contingency of the historical evolution of mathematics.
\footnote{For a discussion of the issue of determinism \emph{vs}
contingency of the historical evolution of \mbox{mathematics}, see
e.g., Mancosu \cite{Ma09} and Bair et al.~\cite[Section~4]{22a}.}

The viability of applying nonstandard analysis (or its axiomatic
versions such as SPOT) to interpreting the procedures of the
historical infinitesimalists depends crucially on the
procedure/foundation distinction.%
\footnote{This point seems to have been overlooked by the authors of
\cite{Ar22b}.  See further in Bascelli et al.~\cite{16a}, Bair et
al.~\cite{21a}.}
%

%this is where the ishiguro stuff is supposed to go

In sum, the theory SPOT has three axioms (in addition to the ZF
axioms): (1) \textbf{S}tandard \textbf{P}art, (2)
n\textbf{O}ntriviality, and (3) \textbf{T}ransfer.  We note the
following.

\begin{enumerate}
\item
\label{ione}
Leibniz emphasized repeatedly that he is working with a relation of
\emph{equality up to} terms that need to be discarded, rather than
with exact equality;
%(see Section~\ref{f29})
he was aware of the need to take the assignable part (in modern
terminology, Standard Part or shadow) of a quantity when appropriate,
as in passing from~$\frac{2x+dx}{a}$ to~$\frac{2x}{a}$ in the
calculation of~$\frac{dy}{dx}$ when~$ay=x^2$; see e.g.,~\cite{Le01c}.
\item
Leibniz postulated inassignable infinitesimals, i.e., Nontriviality.
\item
Leibniz had a Law of Continuity one of whose formulations was that the
rules of the finite succeed in the infinite and vice versa, which
Robinson described as being remarkably close to Transfer.%
\footnote{Robinson \cite[p.\;266]{Ro66}.  See also Section~\ref{s62}.}
Where Leibniz spoke of assignable and inassignable quantities,
Robinson spoke of standard and nonstandard numbers.
\end{enumerate}

Arguably, Berkeley's logical criticism applies to Leibniz's calculus
if and only if it applies to SPOT; see further in
Section~\ref{s9}.%
\footnote{See also \cite[Section~1.1]{22a}.}

\subsection{A note on well-foundedness}
\label{s54}

Since the axioms of SPOT include the \mbox{axioms} of ZF (for
$\in$-formulas), the definitions of~$\N$ and~$\R$ carry over from ZF.
For \mbox{example},
$\N$ could be defined as the smallest inductive set, and~$\R$ by
Dedekind cuts or Cauchy \mbox{sequences}.  There is therefore no
reason to use new symbols for the natural numbers or the reals.  The
same remark applies to all standard concepts (ordinals,
\mbox{cardinals}, etc.).  Some standard concepts (limit, derivative,
\ldots) also have \mbox{nonstandard} \mbox{definitions} in terms of
infinitesimals, which may have their advantages (e.g., in teaching).

The axiomatic approach to NSA does not proceed via models; it
axiomatizes its universe of discourse -- say, by the axioms of SPOT.
The universe of \mbox{discourse} of SPOT satisfies the axiom of
foundation (i.e., every nonempty \emph{set} has an
\mbox{$\in$-minimal} element).  The additional predicate~$\st$ enables
descriptions of certain \emph{non-sets} (classes) that the traditional
mathematician ignores (such as the class of all nonstandard integers).
The universe is not well-founded with respect to these non-sets.

Traditional Platonist mathematicians may find such a viewpoint
difficult to \mbox{accept} because on their view every collection of
elements of some set, no \mbox{matter} how described, should be a set.
But the approach fits well with Hamkins's \mbox{multiverse} view of
set theory.  For example, Fletcher et al.%
\footnote{Fletcher et al.~\cite[Section 7.3, pp.\;228--230]{17f}.}
discuss the relationship between nonstandard analysis and the
Gitman--Hamkins ``toy model'' \cite{Gi10} of the set-theoretic
multiverse.%
\footnote{Hamkins \cite{Ha11}, \cite{Ha12}.}

\section
{From infinite sets to ringinals and recent work}
\label{s9}

In this section, we will provide additional historical context for the
idea of \mbox{formalizing} the procedures of the historical
infinitesimal calculus in the theory SPOT summarized in
Section~\ref{s52}.

\subsection{Positing infinite sets}
\label{s71}
\label{s72}
\label{s73}

As is well known, the existence of an infinite set is an axiom of ZF
as formalized in the 20th century by Zermelo.  Positing infinite sets
was a departure from an ancient philosophical tradition at least as
old as Aristotle (see \cite{Ug22}) that held that only potential
infinity is conceivable.%
\footnote{See further in Bair et al.~\cite[Section~1.3]{22a}.}
Thus, in \emph{Physica} one finds:
\begin{enumerate}\item[]
It is reasonable that there should not be held to be an infinite in
respect of addition such as to surpass every magnitude, but that there
should be thought to be such an infinite in the direction of
division. {\ldots}~It is \mbox{reasonable} too to suppose that in
number there is a limit in the direction of the \mbox{minimum}, and
that in the other direction every multitude is always surpassed.  In
magnitude, on the contrary, every magnitude is surpassed in the
direction of smallness, while in the other direction there is no
infinite magnitude. {\ldots}~in the direction of largeness it is
always possible to think of a larger number; for the number of times a
magnitude can be bisected is infinite. Hence this infinite is
potential, never actual: the number of bisections that can be taken
always surpasses any definite multitude.  But this number is not
separable, and its infinity does not persist but consists in a process
of coming to be, like time and the number of time.  (\emph{Physica}
III 207b3--15)%
\footnote{See e.g.,
\url{https://classics.mit.edu/Aristotle/physics.3.iii.html}.}
\end{enumerate}

Infinite wholes in mathematics were outside the conceptual resources
of both \mbox{Leibniz} and Berkeley; both of them subscribed to the
idea that infinite wholes are contradictory.  Specifically, Leibniz's
analysis of Galileo's paradoxes led him to believe that infinite
wholes contradict Euclid's part-whole principle.%
\footnote{See e.g., Bair et al.~\cite[Section~2.8]{22a}.}

In the Leibnizian calculus, one works with both infinitesimal and
infinite \mbox{numbers}.  Leibniz had a special term for the inverse
of an infinitesimal.  He referred to it as \emph{infinitum terminatum}
(literally: bounded infinity).  The latter was contrasted with
\emph{infinitum \mbox{interminatum}} (unbounded infinity, such as an
unbounded infinite line).%
\footnote{See Bair et al.~\cite[Section 2.2]{21a}.}
\emph{Infinita terminata} were useful in geometry and the calculus,
whereas \emph{infinita \mbox{interminata}} were useless \mbox{because}
contradictory (as mentioned above).  Leibniz's bounded infinity
operates as a kind of ringinal (see Section~\ref{s42}) in the
procedures of the Leibnizian calculus.

\subsection{Postulating predicates and infinite sets}

As noted in Section~\ref{s71}, the modern framework ZF postulates that
an infinite set exists.  Somewhat analogously, one introduces a new
predicate~{\st} and postulating its basic properties, in particular
distinguishing between standard and nonstandard numbers, in modern
axiomatic theories such as IST, HST, or SPOT (see
Section~\ref{s19}).\, Both the axiom of infinity and modern
infinitesimals are 20th century innovations.  Arguably, the acceptance
of infinite sets is neither more nor less justified than the
acceptance of a richer background language incorporating the one-place
relation~{\st} (in addition to the two-place relation~$\in$).
Infinite sets and the properties of the predicate {\st} are similar in
that both are \emph{postulated} rather than \emph{constructed}.
Bed\"urftig and Murawski describe the somewhat paradoxical situation
today as follows:
\begin{enumerate}\item[]
It is astonishing that today everybody believes in the Axiom of
Infinity \mbox{generating} infinitely large sets and infinite cardinal
numbers but \mbox{infinitely} small quantities and infinitesimal
numbers are sometimes indignantly \mbox{rejected}.%
\footnote{Bed\"urftig and Murawski \cite[p.\;201]{Be18}.}
\end{enumerate}
Unlike infinite sets (which have almost no antecedent in mathematics
before \mbox{Cantor}), the predicate `{\st}' has a historical
antecedent in Leibniz's distinction between \mbox{assignable} and
inassignable quantities (see Section~\ref{s13}).

\subsection{Leibniz on violation of Euclid V, Definition 4}
\label{s75}

Leibniz makes it clear that~$dx$ is an element of a non-Archimedean
number system in a 1695 article in \emph{Act.\;Erud.}%
\footnote{Leibniz \cite[p.\;322]{Le95b}.}
in response to Nieuwentijt.  Here Leibniz makes it clear that his
\mbox{incomparable} infinitesimals violate the comparability notion
put forward by \mbox{Euclid} in Book\;V Definition 4, closely related
to the Archimedean property.  Similar \mbox{comments} appeared in a
letter to l'H\^opital%
\footnote{Leibniz \cite[p.\;288]{Le95a}.}
the same year.%
 \footnote{\label{Sc1}It is therefore difficult to agree with Gert
 Schubring's claim that ``{}`the founders of the \mbox{calculus'} had
 not at all created it with a non-Archimedean continuum in mind''
 \cite[p.\;6]{Sc22}, or that ``\mbox{Leibniz} always refused to be
 identified with a foundation based on---rather vaguely
 conceived---infinitesimals'' \cite[p.\,1]{Sc22}.
%(see also   note~\ref{f8}).  
Contrary to Schubring's claim that our group ``for several years has
been leading a \emph{crusade} against the historiography of
mathematics, in particular in its old-fashioned forms'' (ibid.;
emphasis added), we are only ``against'' sloppy historical writing of
the kind \mbox{unfortunately} found in Schubring's publications; see
e.g., B{\l}aszczyk et al.~\cite{17e} and Katz et al.~\cite[notes 46,
  50]{24c}.  Schubring's critique is analyzed in Bair et
al.~\cite[Sections~2--3]{22a}.}
Notes Malet:
\begin{enumerate}
\item[] By the late 17th century {\ldots}~the dividing line between
  numbers and \mbox{continuous} magnitudes was largely gone.%
\footnote{Malet \cite[p.\;213]{Ma12}.}
\end{enumerate}
For details see \cite[Section~3.2]{18a}.  Interpretation of the
Leibnizian calculus is an area of lively debate; see e.g., Eklund
\cite{Ek20} (2020), Esquisabel and Raffo Quintana \cite{Es21} (2021),
Katz et al.~\cite{22c} (2022), Archibald et al.~\cite{Ar22b} (2022),
Bair et al.~\cite{23a} (2023), Katz, Sherry, and Ugaglia~\cite{23h}
(2023), Katz and Kuhlemann \cite{23f} (2023), Katz, Kuhlemann and
Sherry \cite{24b} (2024), Katz, Kuhlemann, Sherry and Ugaglia
\cite{24c} (2024), Ugaglia and Katz \cite{24d} (2024),
Knobloch~\cite{Kn24} (2024), Arthur and Rabouin \cite{Ar24} (2024),
\cite{Ar24b} (2024), Katz \cite{25a} (2025), Kuhlemann \cite{25d}
(2025).

\subsection{Which~$\R$ provides a better background?}

Can one exploit the predicate ``{\st}'' before providing a
construction of a non-Archimedean number system with standard and
nonstandard elements?  The key question here is what is meant by
\emph{construction}.  Leibniz obviously did not possess a modern
set-theoretic construction of \mbox{either} a Weierstrassian real
line~$\R_{W}$, or a Hrbacek--Nelson type real line~$\R_{H\!N}$ (see
\mbox{Section}~\ref{s19}).  Rather, Leibniz presented a collection of
coherent mathematical \mbox{procedures} for \mbox{infinitesimal}
calculus (exploiting the distinction between assignable and
inassignable numbers), and one can ask: comparing theories ZF and SPOT
(see Section~\ref{s52}), which one provides a better background for
formalizing the kind of procedures Leibniz was working with?

%From Tridentine anathemas to Connesian chimeras, infinitesimals
%haven't ceased to fire the imagination.

Following the creation of modern infinitesimal analysis by Robinson in
the 1960s, the field has developed into a vast research area that is
not easy to survey \mbox{meaningfully} in an article on the history of
philosophy.  Many important applications are \mbox{presented} in the
2015 monograph edited by Loeb and Wolf \cite{Loeb2}.  Other recent
exciting work includes Sanders \cite{18l} (2018), Dinis and van den
Berg \cite{Di19} (2019), Goldbring \cite{Go22} (2022), Jin~\cite{Ji23}
(2023), Hrbacek \cite{24a} (2024).

\section{Conclusion}

Due to their apparently elusive nature, infinitesimals as practiced by
Leibniz and others have historically been the subject of reservations
ranging from criticism to outright claims of inherent contradiction.
Some historians have seen modern theories of infinitesimals as
developed by Robinson and others as a \mbox{vindication} of the work
of classical infinitesimalists.  Other historians have argued that
modern \mbox{infinitesimals} exploit mathematical resources undreamt
of by historical \mbox{infinitesimalists}.

The latter criticism is mitigated by means of the distinction,
elaborated e.g., in \cite{17d} and \cite{21a}, between
\emph{foundations} and \emph{procedures}.  Namely, while the
\emph{\mbox{foundational}} aspects of modern infinitesimal frameworks
may involve model theory, \mbox{ultraproducts}, axiom of choice, and
other modern resources, the \emph{procedures} of such theories (such
as the relation of infinite proximity, use of ringinals and associated
infinitesimal \mbox{partitions}) provide a closer fit for the
procedures of the classical \mbox{infinitesimalists} than do the
procedures of the traditional (non-infinitesimal) Weierstrassian
\mbox{framework}.

A more far-reaching response is to exploit axiomatic approaches
outlined in \mbox{Section}~\ref{s19}.  A recent axiomatisation makes
no use of the modern set-theoretic \mbox{resources} mentioned above; a
conservativity result shows that it is possible to avoid any reliance
on the axiom of choice or ultraproducts.  The axioms involved find
close parallels in heuristic principles explicit in Leibniz's
\emph{oeuvre}, such as his dichotomy of assignable and inassignable
number and his Law of Continuity (see Sections~\ref{s13}
and~\ref{s62}).

Nonstandard analysis is sometimes criticized for not being
sufficiently \emph{effective} or for not being a satisfactory model
for real world phenomena, due to its alleged reliance on the axiom of
choice, whereas the elements of the so-called standard models~$\N$ and
$\R$ (constructible in ZF) seem to possess direct referents.  Such
\mbox{objections} have now become obsolete due to the existence of
theories that are conservative extensions of ZF where these standard
models also contain nonstandard/inassignable numbers (like unlimited
numbers and infinitesimals), seen as \emph{fictional} by Leibniz.

The existence of such modern formalisations of the Leibnizian calculus
suggests that infinitesimals are `safer' epistemologically than is
sometimes thought, and \mbox{reveals} the hidden potential of the
notion of the continuum sometimes thought to have been definitively
captured by the mathematical developments from a century and a half
ago.

\section*{Acknowledgments} 
The authors are grateful to John T. Baldwin, Karel Hrbacek, and Jamie
Tappenden for helpful comments.

\end{document}